  \documentclass[onecolumn]{article}
  
  \usepackage{mathtools,amsfonts,amsthm}
  \usepackage[inline]{enumitem}
  \usepackage{cite}
  \usepackage[margin=3cm]{geometry}
  \usepackage[title,titletoc,toc]{appendix}
  
  \def\Rbb{\mathbb{R}} 
  \def\d{\mathrm{d}}
  \def\xbo{\boldsymbol{x}}
  \def\etabo{\boldsymbol{\eta}}
  \def\th#1{${#1}^{th}$}
  \def\smallo#1{\ensuremath{\mathop{}\mathopen{}o\mathopen{}\left(#1\right)}}

  \newtheorem{theoremint}{Theorem}
  \newtheorem{definitionint}{Definition}
  \newtheorem{exampleint}{Example}
  \newtheorem{notationint}{Notation}
  \newtheorem{propositionint}{Proposition}
  \newtheorem{corollaryint}{Corollary}
  \newtheorem{propertyint}{Property}
  \newtheorem{lemmaint}{Lemma}
  \newtheorem*{assumptionint}{Modelling assumptions}
  \newtheorem*{approximationint}{Filtering approximations} 
  
      \newcommand*{\hdelimiter}{\vspace{2pt}}

  \newenvironment{theorem}{\hdelimiter\begin{theoremint}}{\hdelimiter\end{theoremint}}
  \newenvironment{definition}{\hdelimiter\begin{definitionint}}{\hdelimiter\end{definitionint}}

  \newenvironment{lemma}{\hdelimiter\begin{lemmaint}}{\hdelimiter\end{lemmaint}}

  \title{A few calculus rules for chain differentials}
  
  \author{D. E. Clark, J. Houssineau, E. D. Delande}
  
  \date{}
  
\begin{document}

  \maketitle

  \begin{abstract}
    This paper summarizes the core definitions and results regarding the chain differential for functions in locally convex topological vector spaces. In addition, it provides a few elementary calculus rules of practical interest, notably for the differentiation of characteristic functionals in various domains of physical science and engineering.  
  \end{abstract}
    
  
  \section{Functional differentiation}
    In this section we discuss two different forms of differential, the G\^{a}teaux differential~\cite{Gateaux} and the chain differential~\cite{bernhard}. The chain differential, which is similar to the epiderivative~\cite{AubinFran}, is adopted since it is possible to determine a chain rule, yet is not as restrictive as the Fr\'echet derivative.
  
    Results are stated for \emph{locally convex topological vector spaces} which include Banach spaces such as Hilbert and Euclidean spaces, e.g., $\mathbb{R}^n$, as well as spaces of test functions for the study of distributions. This type of space is therefore sufficiently general for most practical applications.

    \subsection{G\^ateaux differential}
      \begin{definition}[G\^ateaux differential]\label{def:Gateaux}
	Let $X$ and $Y$ be locally convex topological vector spaces, and let $\Omega$ be an open subset of $X$ and let $f:\Omega\rightarrow Y$. The G\^ateaux differential at $x\in \Omega$ in the direction $\eta\in X$ is
	\begin{equation}
	  \delta f(x;\eta) := \lim_{\epsilon\to 0} \frac{1}{\epsilon}\big(f(x+\epsilon \eta) - f(x)\big)
	\end{equation}
	when the limit exists. If $\delta f(x;\eta)$ exists for all $\eta\in X$ then $f$ is G\^ateaux differentiable at $x$. The G\^ateaux differential is homogeneous of degree one in $\eta$, so that for all real numbers $\alpha$, $\delta f(x;\alpha \eta) = \alpha\delta f(x;\eta)$.
      \end{definition}

    In Definition \ref{def:Gateaux}, the space $X$ might be a function space. In this case, functions on $X$ can be referred to as \emph{functionals}.

    
    \subsection{Chain differential}
      Due to the lack of continuity properties of the G\^ateaux differential, further constraints are required in order to derive a chain rule. Bernhard \cite{bernhard} proposed a new form of G\^ateaux differential defined with sequences, which he called the chain differential. It is not as restrictive as the Fr\'echet derivative though it is still possible to find a chain rule that maintains the general structure.

      \begin{definition}[Chain differential] \label{def:chain_differential}
	The function $f:X\rightarrow Y$, where $X$ and $Y$ are locally convex topological vector spaces, has a \emph{chain differential} $\delta f(x;\eta)$ at point $x \in X$ in the direction $\eta \in X$ if, for any sequence $\eta_m\rightarrow\eta\in X$, and any sequence of real numbers $\theta_m\rightarrow 0$, it holds that the following limit exists
	\begin{equation} \label{eq:chain_differential}
	  \delta f(x;\eta) := \lim_{m\rightarrow \infty} \dfrac{1}{\theta_m} \big(f(x+\theta_m\eta_m)-f(x)\big).
	\end{equation}
	If $X = X_1\times\ldots\times X_n$, where $\{X_i\}_{i=1}^n$ are locally convex topological vector spaces, $\xbo := (x_1,\ldots,x_n) \in X$, and $\etabo := (\eta_1,\ldots,\eta_n) \in X$, the chain differential $\delta f(\xbo;\etabo)$, if it exists, is also called the \emph{total chain differential} of $f$ at point $\xbo$ in the direction $\etabo$.
      \end{definition}
      
      \begin{definition}[\th{n}-order chain differential]
	The \emph{\th{n}-order chain differential} of $f$ at point $x$ in the sequence of directions $(\eta_i)_{i=1}^n$, is defined recursively with
	\begin{equation} \label{eq:nth_order_chain_differential}
	  \delta^n f\big(x; (\eta_i)_{i=1}^n \big) := \delta \big(y \mapsto \delta^{n-1} f(y; (\eta_i)_{i=1}^{n-1})\big)(x; \eta_n).
	\end{equation}
      \end{definition}
      
      For the sake of simplicity, when there is no ambiguity on the point at which the chain differential is evaluated, the chain differential $\delta f(x;\eta)$ may also be written as $\delta (f(x);\eta)$. The \th{n}-order chain differential \eqref{eq:nth_order_chain_differential} then takes the more compact form
      \begin{equation}
	\delta^n f\big(x; (\eta_i)_{i=1}^n \big) := \delta \big(\delta^{n-1} f(x; (\eta_i)_{i=1}^{n-1}); \eta_n\big).
      \end{equation}
      
      Similarly to the notion of partial derivatives, the notion of chain differential can be defined for appropriate multivariate functions.      
      \begin{definition}[Partial chain differential] \label{def:chain_differential_partial}
	Let $\{X_i\}_{i=1}^n$ and $Y$ be locally convex topological vector spaces. The function $f:X_1\times\ldots\times X_n \rightarrow Y$ has a \emph{partial chain differential} $\delta_i f(\xbo;\eta)$ with respect to the \th{i} variable, at point $\xbo = (x_1,\ldots,x_n) \in X_1\times\ldots\times X_n$ in the direction $\eta \in X_i$ if, for any sequence $\eta_m\rightarrow\eta\in X_i$, and any sequence of real numbers $\theta_m\rightarrow 0$, it holds that the following limit exists
	\begin{equation}
	  \delta_i f(\xbo;\eta) := \lim_{m\rightarrow \infty} \dfrac{1}{\theta_m} \big(f(x_1,\ldots,x_i+\theta_m\eta_m,\ldots,x_n)-f(\xbo)\big).
	\end{equation}
      \end{definition}
      
  \section{Core calculus rules for the chain differential}
    This section summarises the core derivation rules for the chain differential. Note that they all have a counterpart defined for the usual derivative.\newline
  
    \begin{lemma}[Chain rule, from \cite{bernhard}, Theorem 1] \label{lem:chain_rule}
      Let $X$, $Y$ and $Z$ be locally convex topological vector spaces, $f: Y \rightarrow Z$ , $g : X\rightarrow Y$ and $g$ and $f$ have chain differentials at $x$ in the direction $\eta$ and at $g(x)$ in the direction $\delta g(x;\eta)$ respectively. Then the composition $f \circ g$ has a chain differential at point $x$ in the direction $\eta$, given by the chain rule
      \begin{equation} \label{eq:chain_rule}
	\delta (f \circ g)(x;\eta) = \delta f\big(g(x); \delta g(x;\eta)\big).
      \end{equation}
    \end{lemma}
    Note that, unlike its counterpart for the usual derivative, the chain differential of a composition $f \circ g$ does \emph{not} reduce to the product of a chain differential of $f$ and a chain differential of $g$. This key difference has important implications on the structure of the general higher-order chain rule in Theorem~\ref{th:general_chain_rule}.\newline

    Similarly to the usual derivatives, the total chain differential of a multivariate function (see Definition~\ref{def:chain_differential}) can be constructed, in certain conditions defined in the theorem below, as a sum involving its partial chain differentials (see Definition~\ref{def:chain_differential_partial}).
    
    \begin{theorem}[Total chain differential, from \cite{Clark_DE_2013_3}, Theorem 1] \label{th:total_chain}
      Let $\{X_i\}_{i=1}^n$ and $Y$ be locally convex topological vector spaces, $f: X_1\times\ldots\times X_n \rightarrow Y$, $\xbo \in X_1\times\ldots\times X_n$ and $\etabo := (\eta_1,\ldots,\eta_n) \in X_1\times\ldots\times X_n$. If, for $1 \leq i \leq n$, it holds that
      \begin{enumerate}
	\item the partial chain differential $\delta_i f$ exists in a neighbourhood $\Omega \subseteq X_1\times\ldots\times X_n$ of $\xbo$, and,
	\item the function $(x, \eta) \mapsto \delta_i f(x_1,\ldots,x_{i-1},x,x_{i+1},\ldots, x_n; \eta)$ is continuous over $\Omega \times X_i$,
      \end{enumerate}
      then $f$ has a total chain differential at point $\xbo$ in the direction $\etabo$, and it is given by
      \begin{equation}
	\delta f(\xbo;\etabo) = \sum_{i=1}^n \delta_i f(\xbo;\eta_i).
      \end{equation}
    \end{theorem}
    The proof is given in appendix in Section~\ref{subsec:proof_total_chain}. This intermediary result is an important component in the construction of the general higher-order chain rule in Theorem~\ref{th:general_chain_rule}.\newline

    \begin{theorem}[General higher-order chain rule, from \cite{Clark_DE_2013_3}, Theorem 2] \label{th:general_chain_rule}
      Let $X$, $Y$ and $Z$ be  locally convex topological vector spaces, and $f:Y\rightarrow Z$. Assume that $g:X\rightarrow Y$ has higher order chain differentials at point $x$ in all the sequences of directions $(\eta_i)_{i \in I}$, $I \subseteq \{1,\ldots,n\}$. Assume additionally that there exists an open subset $\Omega \subseteq Y$ such that $g(x) \in \Omega$ and $\delta^{|I|} g(x;(\eta_i)_{i \in I}) \in \Omega$, $I \subseteq \{1,\ldots,n\}$. If, for every point $y \in \Omega$ and every sequence $(\xi_i)_{i=1}^m \in \Omega^m$, $1 \leq m \leq n$, it holds that
      \begin{enumerate}
	\item $f$ has a \th{m}-order chain differential at point $y$ in the sequence of directions $(\xi_i)_{i=1}^m$, and,
	\item the conditions of Theorem~\ref{th:total_chain} hold for the function $(y,\xi_1,\ldots,\xi_m) \mapsto \delta^m f(y; (\xi_i)_{i=1}^m)$,
	\item the functions $\xi_j \mapsto \delta^m f(y; (\xi_i)_{i=1}^m)$, $1 \leq j \leq m$, are linear and continuous on $\Omega$,
      \end{enumerate}
      then the \th{n}-order chain differential of the composition $f\circ g$ at point $x$ in the sequence of directions $(\eta_i)_{i=1}^n$ is given by
      \begin{equation}
	\delta^n (f \circ g)(x; (\eta_i)_{i=1}^n) = \sum_{\pi\in \Pi_n} \delta^{|\pi|}f \Big(g(x); \big(\delta^{|\omega|} g(x; (\eta_i)_{i \in \omega})\big)_{\omega \in \pi} \Big),
      \end{equation}
      where $\Pi_n := \Pi(\{1,\ldots,n\})$ denotes the set of the partitions of the index set $\{1,\ldots,n\}$, and $|\pi|$ denotes the cardinality of the set $\pi$.
    \end{theorem}
    The proof is given in appendix in Section~\ref{subsec:proof_general_chain_rule}. The counterpart for usual derivatives is known as F\`{a}a di Bruno's rule \cite{faa2}.\newline
    
    \begin{theorem}[General higher-order product rule] \label{th:product_rule}
      Let $X$ and $Y$ be locally convex topological vector spaces and let $g:X\rightarrow Y$ and $z:X\rightarrow Y$. Assuming that $f$ and $g$ have higher order chain differentials at point $x$ in all the sequences of directions $(\eta_i)_{i \in I}$, $I \subseteq \{1,\ldots,n\}$, then the product $f \cdot g$ has a \th{n}-order chain differential at point $x$ in the sequence of directions $(\eta_i)_{i=1}^n$, and it is given by
      \begin{equation} \label{eq:product_rule}
	\delta^n (f \cdot g)(x; (\eta_i)_{i=1}^n) = \sum_{\pi \subseteq  \{1,\ldots,n\}} \delta^{|\pi|} f(x; (\eta_i)_{i \in \pi}) \delta^{n-|\pi|} g(x; (\eta_i)_{i \in \pi^c})
      \end{equation}
      where $\pi^c := \{1, \ldots, n\} \setminus \pi$ denotes the complement of $\pi$ in $\{1, \ldots, n\}$.
    \end{theorem}
    The proof is given in appendix in Section~\ref{subsec:proof_product_rule}. The counterpart for usual derivatives is known as Leibniz' rule.
    
  \section{Practical derivations for the chain rule}
    This section provides specific applications of the chain rule given in Lemma~\ref{lem:chain_rule} in which the outer function of the composition $f \circ g$ assumes a specific form, commonly encountered in practical derivations. 
  
    \begin{theorem}[Practical derivations for the chain rule] \label{th:chain_rule_usual}
      Let $X$, $Y$ and $Z$ be locally convex topological vector spaces, $f : Y \rightarrow Z$ , $g : X\rightarrow Y$. Assume additionally that $g$ has a chain differential at some point $x$ in some direction $\eta$. Then:\newline
      
      a) if $f$ is a continuous linear function $\ell$, then the composition $f \circ g$ has a chain differential at point $x$ in the direction $\eta$, and it is given by
	\begin{equation}
	  \delta (\ell \circ g)(x;\eta) = \ell\big(\delta g(x;\eta)\big).
	\end{equation}Leibniz' rule
	
      b) if $f$ is the \th{k} power function $y \mapsto y^k$, $k > 0$, then the composition $f \circ g$ has a chain differential at point $x$ in the direction $\eta$, and it is given by
      \begin{equation}
	\delta \big((y \mapsto y^k) \circ g\big)(x;\eta) = k\big(g(x)\big)^{k-1}\delta g(x;\eta).
      \end{equation}
	
      c) if $f$ is the exponential function $\exp$, then the composition $f \circ g$ has a chain differential at point $x$ in the direction $\eta$, and it is given by
      \begin{equation}
	\delta (\exp \circ g)(x;\eta) = \exp\big(g(x)\big)\delta g(x;\eta).
      \end{equation}
      
    \end{theorem}
    The proof is given in appendix in Section~\ref{subsec:proof_chain_rule_usual}.

  \begin{appendices}
  \section{Proofs} \label{sec:proof}
    \subsection{Total chain differential (Theorem \ref{th:total_chain})} \label{subsec:proof_total_chain}
      \begin{proof}
	The result is proved in the case $n=2$ from which the general case can be straightforwardly deduced. Let us fix a point $\xbo = (x,y) \in X_1 \times X_2$ and a direction $\etabo = (\eta,\xi) \in X_1 \times X_2$, such that $f$ has partial chain differentials $\delta_1 f(\xbo;\eta)$ and $\delta_2 f(\xbo;\xi)$. Let us then fix arbitrary sequences of directions $\eta_m \rightarrow \eta \in X_1$, $\xi_m \rightarrow \xi \in X_2$, and an arbitrary sequence of real numbers $\theta_m \rightarrow 0$. For $m \leq 0$ we can write
	\begin{subequations}
	  \begin{align}
	    \theta^{-1}_m\big[f(\xbo + \theta_m\etabo) - f(\xbo)\big] &= \theta^{-1}_m\big[f(x+\theta_m \eta_m,y + \theta_m \xi_m) - f(x,y)\big]
	    \\
	    &= \theta^{-1}_m\big[g_1(y + \theta_m\xi_m) - g_1(y) \big] + \theta^{-1}_m\big[g_2(x+\theta_m \eta_m) - g_2(x)\big], \label{eq:decomposition_1}
	  \end{align}
	\end{subequations}
	where we define $g_1(y)$ and $g_2(x)$ as follows:
	\begin{equation}
	  \begin{dcases}
	    g_1(y) &= f(x+\theta_m \eta_m,y),
	    \\
	    g_2(x) &= f(x,y).
	  \end{dcases}
	\end{equation}
	Given $\theta_m \neq 0$, define $h:\Rbb \rightarrow \Rbb$ as $h(t) = g_1(y + t\xi_m)$. From the mean value theorem for real-valued functions, there exists $c_y \in [0,\theta_m]$ such that
	\begin{subequations}
	  \begin{align}
	    \theta_m^{-1}\big[h(\theta_m) - h(0)\big] &= \left.\dfrac{\d h}{\d t}\right|_{t = c_y}
	    \\
	    &= \delta h(c_y;1),
	  \end{align}
	\end{subequations}
	which, when replacing $h(t)$ by $g_1(y + t\xi_m)$, can be rewritten
	\begin{subequations}
	  \begin{align}
	    \theta_m^{-1}\big[g_1(y + \theta_m\xi_m) - g_1(y)\big] &= \delta \big(g_1(y + c_y\xi_m); 1\big)
	    \\
	    &= \delta g_1(y + c_y \xi_m; \xi_m), \label{eq:decomposition_2a}
	  \end{align}
	\end{subequations}
	where Lemma~\ref{lem:chain_rule} has been used to obtain the last equality. Similarly for $g_2(x)$, there exists $c_x \in [0,\theta_m]$ such that
	\begin{equation}
	  \theta^{-1}_m\left[g_2(x+\theta_m \eta_m) - g_2(x)\right] = \delta g_2(x+c_x \eta_m; \eta_m). \label{eq:decomposition_2b}
	\end{equation}
	Let us now prove that the limit of the term
	\begin{equation} \label{eq:int_step}
	  \Big|\theta^{-1}_m\big[f(\xbo + \theta_m\etabo) - f(\xbo)\big] - \delta_1 f(\xbo; \eta_m) - \delta_2 f(\xbo; \xi_m) \Big|
	\end{equation}
	is equal to $0$ when $r\rightarrow\infty$. Substituting \eqref{eq:decomposition_2a} and \eqref{eq:decomposition_2b} into \eqref{eq:decomposition_1}, \eqref{eq:int_step} becomes
	\begin{equation} \label{eq:last_step}
	  \big|\delta g_2(x + c_x \eta_m;\eta_m) + \delta g_1(y + c_y \xi_m;\xi_m) - \delta_1 f(\xbo; \eta_m) - \delta_2 f(\xbo; \xi_m) \big|.
	\end{equation}
	By the triangle inequality, \eqref{eq:last_step} is bounded above by the following summation
	\begin{equation} \label{eq:bound}
	  \big|\delta g_2(x+c_x \eta_m;\eta_m) - \delta_1 f(\xbo; \eta_m) \big| + \big| \delta g_1(y + c_y \xi_m;\xi_m) - \delta_2 f(\xbo; \xi_m) \big|.
	\end{equation}
	Substituting $g_1$ and $g_2$ with $f$, the bound \eqref{eq:bound} becomes
	\begin{align}
	  \big|\delta_1 f(x+c_x \eta_m, y; \eta_m) - \delta_1 f(\xbo; \eta_m)\big| + \big|\delta_2 f(x, y + c_y \xi_m; \xi_m) - \delta_2 f(\xbo; \xi_m)\big|,
	\end{align}
	which tends to $0$ when $m\rightarrow\infty$ because of the continuity of the functions $(z, \nu) \mapsto \delta_1 f(z,y;\nu)$ and $(z, \nu) \mapsto \delta_2 f(x,z;\nu)$. Thus, it holds that
	\begin{equation}
	  \lim_{m \rightarrow \infty} \Big|\theta^{-1}_m\big[f(\xbo + \theta_m\etabo) - f(\xbo)\big] - \delta_1 f(\xbo; \eta_m) - \delta_2 f(\xbo; \xi_m) \Big| = 0,
	\end{equation}
	that is, $f$ has a total chain differential in point $\xbo$ in direction $\etabo$, and it is such that
	\begin{equation}
	  \delta f(\xbo; \etabo) = \delta_1 f(\xbo; \eta) + \delta_2 f(\xbo; \xi),
	\end{equation}
	which is equivalent to the Proposition $3$ in \cite{bernhard}.
      \end{proof}
      
    \subsection{General higher-order chain rule (Theorem \ref{th:general_chain_rule})} \label{subsec:proof_general_chain_rule}
      \begin{proof}
	The proof is constructed by induction on the number of directions $n$. Lemma \ref{lem:chain_rule} gives the base case $n = 1$. For the induction step, we apply the differential operator to the case $n$ to give the case $n+1$ and show that it involves a summation over partitions of elements $\eta_1,\ldots,\eta_{n+1}$ in the following way
	\begin{equation} \label{eq:faaProof1}
	  \delta^{n+1} (f\circ g)(x; (\eta_i)_{i=1}^{n+1}) = 
	  \sum_{\pi\in \Pi_n} \delta \bigg(u \mapsto \delta^{|\pi|}f \Big(g(u); \big(\delta^{|\omega|} g(u; (\eta_j)_{j \in \omega})\big)_{\omega \in \pi} \Big)\bigg)(x; \eta_{n+1}).
	\end{equation}
	
	The main objective in this proof is to calculate a term of the summation on the right-hand side of \eqref{eq:faaProof1}, of the form
	\begin{equation} \label{eq:faaProof2}
	  \delta \Big(u \mapsto \delta^k f \big(g(u); (h_i(u))_{i=1}^k \big)\Big)(x; \eta).
	\end{equation}

	The additional differentiation with respect to $\eta$ applies to every function on $X$, i.e. to the function $g$ and to the functions $h_i$, $1\leq i\leq k$. To highlight the structure of this result, we can define a multi-variate function $F$ such that
	\begin{equation}
	  \begin{aligned}
	    F:\quad Y^{k+1} &\to Z
	    \\
	    (y_0,\ldots,y_k) &\mapsto \delta^k f(y_0; (y_i)_{i=1}^k)
	  \end{aligned},
	\end{equation}
	so that \eqref{eq:faaProof2} can be rewritten as $\delta \bigg(F \circ \Big(u \mapsto \big(g(u),h_1(u),\ldots,h_k(u)\big)\Big)\bigg)(x;\eta)$, which is equal to
	\begin{equation}
	  \delta F\bigg(g(x),h_1(x),\ldots,h_k(x); \delta\Big(u \mapsto \big(g(u),h_1(u),\ldots,h_k(u)\big)\Big)(x; \eta)\bigg), \label{eq:dfb_decomp_1}
	\end{equation}
	using Lemma~\ref{lem:chain_rule}. Let us focus on the direction $\delta\Big(u \mapsto \big(g(u),h_1(u),\ldots,h_k(u)\big)\Big)(x; \eta)$ in \eqref{eq:dfb_decomp_1}. Let us define a sequence $\eta_m \rightarrow \eta \in X$ and a sequence of real numbers $\theta_m \rightarrow 0$. Using the definition of the chain differential \eqref{eq:chain_differential}, we can write
	\begin{subequations}
	  \begin{align}
	    &\delta\Big(u \mapsto \big(g(u),h_1(u),\ldots,h_k(u)\big)\Big)(x; \eta) \nonumber
	    \\
	    &= \lim_{m \rightarrow \infty} \theta_m^{-1}\Big[\big(g(x + \theta_m\eta_m),h_1(x + \theta_m\eta_m),\ldots,h_k(x + \theta_m\eta_m)\big) - \big(g(x),h_1(x),\ldots,h_k(x)\big)\Big]
	    \\
	    &= \big(\delta g(x;\eta), \delta h_1(x;\eta), \ldots, \delta h_k(x;\eta)\big), \label{eq:dfb_decomp_2}
	  \end{align}
	\end{subequations}
	where the last equality is given by the definition of the chain differential \eqref{eq:chain_differential} applied to the functions $g$ and $h_i$, $1 \leq i \leq k$. Substituting \eqref{eq:dfb_decomp_2} into \eqref{eq:dfb_decomp_1} and applying Theorem~\ref{th:total_chain}, \eqref{eq:faaProof2} becomes
	\begin{equation}
	  \delta_1 F \big(g(x),h_1(x),\ldots,h_k(x); \delta g(x;\eta)\big) + \sum_{i=1}^k \delta_{i+1} F \big(g(x),h_1(x),\ldots,h_k(x); \delta h_i(x;\eta)\big). \label{eq:dfb_decomp_3}
	\end{equation}
	
	\begin{itemize}
	  \item Consider the first term of the summation in \eqref{eq:dfb_decomp_3}:
	  \begin{equation}
	    \delta_1 F \big(g(x),h_1(x),\ldots,h_k(x); \delta g(x;\eta)\big).
	  \end{equation}
	  Using the definition of $F$, it can be written as $\delta \Big(y \mapsto \delta^k f\big(y; (h_i(x))_{i=1}^k\big)\Big) \big(g(x);\delta g(x;\eta)\big)$, which is equal to
	  \begin{equation} \label{eq:dfb_decomp_4}
	    \delta^{k+1} f \big(g(x); h_1(x),\ldots,h_k(x),\delta g(x;\eta)\big),
	  \end{equation}
	  by definition of the \th{(k+1)}-order chain differential.
	  
	  \item Now consider any other term in \eqref{eq:dfb_decomp_3}:
	  \begin{equation} \label{eq:otherTerm}
	     \delta_{i+1} F \big(g(x),h_1(x),\ldots,h_k(x); \delta h_i(x;\eta)\big).
	  \end{equation}
	  Using the definition of $F$, it can be written as
	  \begin{equation} \label{eq:dfb_decomp_5}
	    \delta \Big(y \mapsto \delta^k f\big(g(x); h_1(x),\ldots,y,\ldots,h_k(x)\big)\Big) \big(h_i(x);\delta h_i(x;\eta)\big).
	  \end{equation}
	  Let us define a sequence $\nu_m \rightarrow \delta h_i(x;\eta) \in Y$ and a sequence of real numbers $\theta_m \rightarrow 0$. Using the definition of the chain differential \eqref{eq:chain_differential}, \eqref{eq:dfb_decomp_5} becomes
	  \begin{subequations}
	    \begin{align}
	      &\lim_{m \rightarrow \infty}\theta_m^{-1}\Big[\delta^k f\big(g(x); h_1(x),\ldots,h_i(x) + \theta_m\nu_m,\ldots,h_1(x)\big) - \delta^k f\big(g(x); h_1(x),\ldots,h_i(x),\ldots,h_1(x)\big)\Big] \nonumber
	      \\
	      &= \lim_{m \rightarrow \infty} \delta^k f\big(g(x); h_1(x),\ldots,\nu_m,\ldots,h_1(x)\big)
	      \\
	      &= \delta^k f\big(g(x); h_1(x),\ldots,\delta h_i(x;\eta),\ldots,h_1(x)\big), \label{eq:dfb_decomp_6}
	    \end{align}
	  \end{subequations}
	  while the first equality exploits the linearity of the function $y \mapsto \delta^k f\big(g(x); h_1(x),\ldots,y,\ldots,h_k(x)\big)$, and the second equality its continuity on $\Omega$.
	\end{itemize}
	Substituting \eqref{eq:dfb_decomp_4} and \eqref{eq:dfb_decomp_6} in \eqref{eq:dfb_decomp_3}, \eqref{eq:faaProof2} becomes
	\begin{equation}
	  \delta^{k+1} f \big(g(x); h_1(x),\ldots,h_k(x),\delta g(x;\eta)\big) + \sum_{i=1}^k \delta^k f\big(g(x); h_1(x),\ldots,\delta h_i(x;\eta),\ldots,h_1(x)\big). \label{eq:dfb_decomp_7}
	\end{equation}

	Considering $\eta := \eta_{n+1}$ and $h_i(x) := \delta^{|\omega_i|} g \big(x; (\eta_j)_{j \in \omega_i}\big)$ and replacing the result \eqref{eq:dfb_decomp_7} into \eqref{eq:faaProof1}, we find
	\begin{subequations}
	  \begin{align}
	    &\delta^{n+1}(f\circ g)(x; (\eta_i)_{i=1}^{n+1}) \nonumber
	    \\
	    &= \sum_{\pi\in \Pi_n} \delta^{|{\pi}|+1} f \Big(g(x); \big(\delta^{|\omega|} g (x; (\eta_j)_{j \in \omega})\big)_{\omega \in \pi \cup \{\{n+1\}\}}\Big) \nonumber
	    \\
	    &\hspace{5cm}+ \sum_{\pi\in \Pi_n} \sum_{\nu \in \pi} \delta^{|{\pi}|} f \Big(g(x); \big(\delta^{|\omega|} g(x; (\eta_j)_{j \in \omega})\big)_{\omega \in \pi \setminus \{\nu\} \cup \{\nu \cup \{n+1\}\}}\Big) \label{eq:partition}
	    \\
	    &= \sum_{\pi\in \Pi_{n+1}} \delta^{|\pi|}f \Big(g(x); \big(\delta^{|\omega|} g(x; (\eta_j)_{j \in \omega})\big)_{\omega \in \pi} \Big).
	  \end{align}
	\end{subequations}

	Following a similar argument used for the recursion of Stirling numbers of the second kind and their relation to Bell numbers~\cite[p74]{RStanleyV1}, the result above can be viewed as a means of generating all partitions of $n+1$ elements from all partitions of $n$ elements: The first term in Eq.~\eqref{eq:partition} corresponds to the creation of a new element to the partition $\pi \in \Pi_n$, containing only $n+1$, and each term in the second summation appends $n+1$ to one of the existing element $\nu$ of the partition $\pi$. This argument follows similar arguments previously used for ordinary and partial derivatives~\cite{ma,huang,hardy}. Hence the result is proved by induction.
      \end{proof}

    \subsection{General higher-order product rule (Theorem \ref{th:product_rule})} \label{subsec:proof_product_rule}
      \begin{proof}
	The proof is constructed by induction on the number of directions $n$.
	
	a) Case $n = 0$.\newline
	We can write immediately
	\begin{subequations}
	  \begin{align}
	    \sum_{\pi \subseteq \emptyset} \delta^{|\pi|} f(x; (\eta_i)_{i \in \pi}) \delta^{0 - |\pi|} g(x; (\eta_i)_{i \in \emptyset \setminus \pi}) &= \delta^0 f(x)\delta^0 g(x)
	    \\
	    &= f(x)g(x)
	    \\
	    &= (f \cdot g)(x)
	    \\
	    &= \delta^0(f \cdot g)(x).
	  \end{align}
	\end{subequations}

	b) Case $n = 1$.\newline
	Let us fix a point $x \in X$ and a direction $\eta_1 \in X$, such that both $f$ and $g$ have a first-order chain differential at point $x$ in direction $\eta_1$. Let us then fix a sequence $\eta_{1, m} \rightarrow \eta_1 \in X$, and a sequence of real numbers $\theta_m \rightarrow 0$. Since $f$ has a first-order chain differential at point $x$, it is continuous in $x$ and thus
	\begin{equation}
	  \lim_{m\rightarrow \infty} f(x + \theta_m\eta_{1, m}) = f(x). \label{eq:condition_1}
	\end{equation}
	Since both $f$ and $g$ have a first-order chain differential at point $x$ in direction $\eta_1$, we have \cite{bernhard}
	\begin{align}
	  \lim_{m\rightarrow \infty} \theta_m^{-1} \big[f(x + \theta_m\eta_{1, m}) - f(x)\big] &= \delta f(x; \eta_1), \label{eq:condition_2}
	  \\
	  \lim_{m\rightarrow \infty} \theta_m^{-1} \big[g(x + \theta_m\eta_{1, m}) - g(x)\big] &= \delta g(x; \eta_1). \label{eq:condition_3}
	\end{align}
	From \eqref{eq:condition_1}, \eqref{eq:condition_2}, and \eqref{eq:condition_3}, it holds that
	\begin{subequations}
	  \begin{align}
	    &\lim_{m\rightarrow \infty} \theta_m^{-1} \big[(f \cdot g)(x + \theta_m\eta_{1,m}) - (f \cdot g)(x)\big] \nonumber
	    \\
	    &= \lim_{m\rightarrow \infty} \theta_m^{-1} \big[f(x + \theta_m\eta_{1,m})g(x + \theta_m\eta_{1,m}) - f(x + \theta_m\eta_{1,m})g(x) + f(x + \theta_m\eta_{1,m})g(x) - f(x)g(x)\big]
	    \\
	    &= \lim_{m\rightarrow \infty} f(x + \theta_m\eta_{1,m})\lim_{m\rightarrow \infty} \theta_m^{-1} \big[g(x + \theta_m\eta_{1,m}) - g(x)\big] \nonumber
	    \\
	    &\hspace{7cm}+ \lim_{m\rightarrow \infty} \theta_m^{-1} \big[f(x + \theta_m\eta_{1,m}) - f(x)\big]g(x),
	    \\
	    &= f(x)\delta g(x; \eta_1) + \delta f(x; \eta_1)g(x).
	  \end{align}
	\end{subequations}
	That is, $f \cdot g$ has a first-order chain differential at point $x$ in direction $\eta_1$ and it is such that
	\begin{equation}
	   \delta (f \cdot g)(x; \eta_1) = f(x)\delta g(x; \eta_1) + \delta f(x; \eta_1)g(x).
	\end{equation}

	c) Case $n \geq 2$.\newline
	Let us fix a point $x \in X$ and a sequence of directions $(\eta_i)_{i = 1}^n \in X^n$ such that both $f$ and $g$ have higher order chain differentials at point $x$ in all sequences of directions $(\eta_i)_{i \in I}$, $I \subseteq \{1,\ldots,n\}$. We can then write
	\begin{subequations}
	  \begin{align}
	    &\sum_{\pi \subseteq \{1,\ldots,n\}} \delta^{|\pi|} f(x; (\eta_i)_{i \in \pi}) \delta^{n - |\pi|} g(x; (\eta_i)_{i \in \{1,\ldots,n\} \setminus \pi}) \nonumber
	    \\
	    &= \sum_{\pi \subseteq \{1,\ldots,n-1\}} \delta^{|\pi| + 1} f(x; (\eta_i)_{i \in \pi \cup \{n\}}) \delta^{n - 1 - |\pi|} g(x; (\eta_i)_{i \in \{1,\ldots,n-1\} \setminus \pi}) \nonumber
	    \\
	    &\hspace{4cm} + \sum_{\pi \subseteq \{1,\ldots,n-1\}} \delta^{|\pi|} f(x; (\eta_i)_{i \in \pi}) \delta^{n - 1 - |\pi| + 1} g(x; (\eta_i)_{i \in \{1,\ldots,n-1\} \setminus \pi \cup \{n\}})
	    \\
	    &= \sum_{\pi \subseteq \{1,\ldots,n-1\}} \delta\big(\delta^{|\pi|} f(x; (\eta_i)_{i \in \pi}); \eta_n\big) \delta^{n - 1 - |\pi|} g(x; (\eta_i)_{i \in \{1,\ldots,n-1\} \setminus \pi}) \nonumber
	    \\
	    &\hspace{4cm} + \sum_{\pi \subseteq \{1,\ldots,n-1\}} \delta^{|\pi|} f(x; (\eta_i)_{i \in \pi}) \delta\big(\delta^{n - 1 - |\pi|} g(x; (\eta_i)_{i \in \{1,\ldots,n-1\} \setminus \pi}); \eta_n\big)
	    \\
	    &= \sum_{\pi \subseteq \{1,\ldots,n-1\}} \delta \big(\delta^{|\pi|} f(x; (\eta_i)_{i \in \pi})\delta^{n - 1 - |\pi|} g(x; (\eta_i)_{i \in \{1,\ldots,n-1\} \setminus \pi}); \eta_n\big)
	    \\
	    &= \delta \Big(\sum_{\pi \subseteq \{1,\ldots,n-1\}} \delta^{|\pi|} f(x; (\eta_i)_{i \in \pi})\delta^{n - 1 - |\pi|} g(x; (\eta_i)_{i \in \{1,\ldots,n-1\} \setminus \pi}); \eta_n\Big)
	    \\
	    &= \delta \big(\delta^{n-1}(f \cdot g)(x; (\eta_i)_{i = 1}^{n-1}); \eta_n\big),
	  \end{align}
	\end{subequations}
	where the last equality was obtained by exploiting case $n-1$. Thus, $f \cdot g$ has a \th{n}-order chain differential at point $x$ in directions $\eta_1,\ldots,\eta_n$ and it is such that
	\begin{equation}
	   \delta^n (f \cdot g)(x; (\eta_i)_{i = 1}^{n}) = \sum_{\pi \subseteq \{1,\ldots,n\}} \delta^{|\pi|} f(x; (\eta_i)_{i \in \pi}) \delta^{n - |\pi|} g(x; (\eta_i)_{i \in \{1,\ldots,n\} \setminus \pi}).
	\end{equation}
	This ends the proof by induction.
      \end{proof}
      
    \subsection{Practical derivations for the chain rule (Theorem \ref{th:chain_rule_usual})} \label{subsec:proof_chain_rule_usual}
      \begin{proof}
	Let us fix a point $x \in X$ and a direction $\eta \in X$, such that $g$ has a chain differential at point $x$ in direction $\eta$. Let us then fix a sequence $\nu_m \rightarrow \delta g(x;\eta) \in Y$, and a sequence of real numbers $\theta_m \rightarrow 0$.\newline
	
	a) Let us assume that $f$ is a continuous linear function $\ell$. For any $m \geq 0$ we can write
	\begin{subequations}
	  \begin{align}
	    \theta_m^{-1}\big[\ell(g(x) + \theta_m\nu_m) - \ell(g(x))\big] &= \theta_m^{-1}\big[\ell(g(x)) + \theta_m\ell(\nu_m) - \ell(g(x))\big]
	    \\
	    &= \ell(\nu_m).
	  \end{align}
	\end{subequations}
	Thus, it holds that
	\begin{subequations}
	  \begin{align}
	    \lim_{m \rightarrow \infty} \theta_m^{-1}\big[\ell(g(x) + \theta_m\nu_m) - \ell(g(x))\big] &= \lim_{m \rightarrow \infty}\ell(\nu_m)
	    \\
	    &= \ell\big(\lim_{m \rightarrow \infty}\nu_m\big)
	    \\
	    &= \ell(\delta g(x;\eta)).
	  \end{align}
	\end{subequations}
	Thus, using the definition of the chain differential \eqref{eq:chain_differential}, we have
	\begin{equation}
	  \delta \ell\big(g(x); \delta g(x;\eta)\big) = \ell\big(\delta g(x;\eta)\big).
	\end{equation}
	Using the chain rule \eqref{eq:chain_rule} ends the proof.\newline
	
	b) Let us assume that $f$ is the \th{k}-power function $y \mapsto y^k$ for some $k > 0$. For any $m \geq 0$ we can write
	\begin{subequations}
	  \begin{align}
	    &\theta_m^{-1}\Big[\big(g(x) + \theta_m\nu_m\big)^k - \big(g(x)\big)^k\Big] \nonumber
	    \\
	    &= \theta_m^{-1}\Big[\sum_{p = 0}^k {k \choose p}\big(g(x)\big)^p(\theta_m\nu_m)^{k-p} - \big(g(x)\big)^k\Big]
	    \\
	    &= \theta_m^{-1}\Big[\big(g(x)\big)^k + k\big(g(x)\big)^{k-1}\theta_m\nu_m + \sum_{p = 0}^{k-2} {k \choose p}\big(g(x)\big)^p(\theta_m\nu_m)^{k-p} - \big(g(x)\big)^k\Big]
	    \\
	    &= k\big(g(x)\big)^{k-1}\nu_m + \sum_{p = 0}^{k-2} {k \choose p}\big(g(x)\big)^p(\theta_m)^{k-p-1}(\nu_m)^{k-p}.
	  \end{align}
	\end{subequations}
	Thus, it holds that
	\begin{subequations}
	  \begin{align}
	    &\lim_{m \rightarrow \infty} \theta_m^{-1}\Big[\big(g(x) + \theta_m\nu_m\big)^k - \big(g(x)\big)^k\Big] \nonumber
	    \\
	    &= k\big(g(x)\big)^{k-1}\underbrace{\lim_{m \rightarrow \infty} \nu_m}_{=\delta g(x;\eta)} + \sum_{p = 0}^{k-2} {k \choose p}\big(g(x)\big)^p\underbrace{\lim_{m \rightarrow \infty} (\theta_m)^{k-p-1}}_{=0}\underbrace{\lim_{m \rightarrow \infty} (\nu_m)^{k-p}}_{= \delta g(x;\eta)^{k-p}}
	    \\
	    &= k\big(g(x)\big)^{k-1}\delta g(x;\eta).
	  \end{align}
	\end{subequations}
	Thus, using the definition of the chain differential \eqref{eq:chain_differential}, we have
	\begin{equation}
	  \delta (y \mapsto y^k)\big(g(x); \delta g(x;\eta)\big) = k\big(g(x)\big)^{k-1}\delta g(x;\eta).
	\end{equation}
	Using the chain rule \eqref{eq:chain_rule} ends the proof.\newline
	
	c) Let us assume that $f$ is the exponential function. For any $m \geq 0$ we can write
	\begin{subequations}
	  \begin{align}
	    \theta_m^{-1}\Big[\exp\big(g(x) + \theta_m\nu_m\big) - \exp\big(g(x)\big)\Big] &= \theta_m^{-1}\exp\big(g(x)\big)\big[\exp(\theta_m\nu_m) - 1\big]
	    \\
	    &= \theta_m^{-1}\exp\big(g(x)\big)\big[\theta_m\nu_m + \smallo{\theta_m\nu_m}\big].
	  \end{align}
	\end{subequations}
	Thus, it holds that
	\begin{subequations}
	  \begin{align}
	    \lim_{m \rightarrow \infty} \theta_m^{-1}\Big[\exp\big(g(x) + \theta_m\nu_m\big) - \exp\big(g(x)\big)\Big] &= \exp\big(g(x)\big)\big[\underbrace{\lim_{m \rightarrow \infty} \nu_m}_{=\delta g(x;\eta)} + \underbrace{\lim_{m \rightarrow \infty} \theta_m^{-1}\smallo{\theta_m\nu_m}}_{=0}\big]
	    \\
	    &= \exp\big(g(x)\big)\delta g(x;\eta).
	  \end{align}
	\end{subequations}
	Thus, using the definition of the chain differential \eqref{eq:chain_differential}, we have
	\begin{equation}
	  \delta \exp\big(g(x); \delta g(x;\eta)\big) = \exp\big(g(x)\big)\delta g(x;\eta).
	\end{equation}
	Using the chain rule \eqref{eq:chain_rule} ends the proof.
      \end{proof}
    \end{appendices}

  \bibliographystyle{plain}
  \bibliography{bibHierSys.bib}

\end{document}